\documentclass[11pt]{article}
\usepackage{amsmath, amsthm}
\usepackage{latexsym}
\usepackage{graphicx}
\makeatletter
\newfont{\footsc}{cmcsc10 at 8truept}
\newfont{\footbf}{cmbx10 at 8truept}
\newfont{\footrm}{cmr10 at 10truept}
\makeatother
\newtheorem{theorem}{\bf Theorem}

\newtheorem{lemma}{\bf Lemma}

\begin{document}
\title{An Inequality for Ratios of Gamma Functions}

\author{Yaming Yu\\
\small Department of Statistics\\[-0.8ex]
\small University of California\\[-0.8ex]
\small Irvine, CA 92697, USA\\[-0.8ex]
\small \texttt{yamingy@uci.edu}}

\date{}
\maketitle

\begin{abstract}
Let $\Gamma(x)$ denote Euler's gamma function.  The following inequality
is proved: for $y>0$ and $x>1$ we have
$$\frac{[\Gamma(x+y+1)/\Gamma(y+1)]^{1/x}}{[\Gamma(x+y+2)/\Gamma(y+1)]^{1/(x+1)}}
< \sqrt{\frac{x+y}{x+y+1}}.$$
The inequality is reversed if $0<x<1$.  This resolves an open problem of
Guo and Qi (2003).

\end{abstract}

\section{Main result}
Recently there has been considerable interest in inequalities concerning Euler's gamma function 
$$\Gamma(z)=\int_0^\infty u^{z-1} e^{-u}\, {\rm d}u,\quad z>0.$$  
This note resolves a problem left open by Guo and Qi \cite{GQ} on this topic.

\begin{theorem}
\label{main}
If $y>0$ and $x>1$, then 
\begin{equation}
\label{eqn1}
\frac{[\Gamma(x+y+1)/\Gamma(y+1)]^{1/x}}{[\Gamma(x+y+2)/\Gamma(y+1)]^{1/(x+1)}} < \sqrt{\frac{x+y}{x+y+1}}.
\end{equation}
If $y>0$ and $0<x<1$, then the inequality is reversed.
\end{theorem}

Note that, by the recursion $\Gamma(z+1)=z\Gamma(z)$, the two sides of (\ref{eqn1}) are equal if $x=1$.  Also, Guo and Qi 
(\cite{GQ}, Theorem 2) gave $(x+y+1)/(x+y+2)$ as a lower bound for the left hand side of (\ref{eqn1}); see Chen and Qi 
\cite{CQ} and the references therein for related work.

\section{Preliminaries}
As usual the digamma function is denoted by $\psi(z)={\rm d}\log \Gamma(z)/{\rm d}z$.  Among the basic properties of 
$\Gamma(z)$ and $\psi(z)$, we shall make use of the recursion
\begin{equation}
\label{psi}
\psi(z+1)=\psi(z)+\frac{1}{z}
\end{equation}
and the asymptotic expansions \cite{AS}
\begin{align}
\label{lgam}
\log\Gamma(z) &= \left(z-\frac{1}{2}\right)\log(z) -z+\frac{1}{2}\log(2\pi)+O(z^{-1}),\\
\nonumber
\psi'(z) &= \frac{1}{z}+\frac{1}{2z^2}+\frac{1}{6z^3}+O(z^{-5}),\\
\nonumber
\psi''(z) &= -\frac{1}{z^2}-\frac{1}{z^3}-\frac{1}{2z^4}+\frac{1}{6z^6}+O(z^{-8}),
\end{align}
valid as $z\rightarrow \infty$.  A result of Alzer (\cite{A}, Theorem 8) implies the following bounds: for $z>0$
\begin{align}
\label{bp}
\psi'(z) &>\frac{1}{z}+\frac{1}{2z^2};\\
\label{bpp}
\psi''(z) &> -\frac{1}{z^2}-\frac{1}{z^3}-\frac{1}{2z^4}.
\end{align}
See also Qi et al.\ \cite{QCCG} and Koumandos (\cite{K}, Theorem 1).  The bounds (\ref{bp}) and (\ref{bpp}) are crucial in our 
proof of Theorem \ref{main}. 

\section{Proof of Theorem \ref{main}}
For $x>0$ and $y>0$ denote
\begin{align*}
f(x, y) &=\frac{\log\Gamma(x+y+1)-\log\Gamma(y+1)}{x}-\frac{1}{2} \log(x+y);\\
g(x, y) &=f(x+1, y)-f(x, y).
\end{align*}
Equivalently we need to show 
\begin{align*}
g(x, y)>0,\quad  &x>1;\\
g(x, y)<0,\quad  &0<x<1. 
\end{align*}
As direct assessment of the sign of $g(x, y)$ appears difficult, let us consider the function $\partial 
g(x, y)/\partial y$ instead.  Using $\psi(x+y+2)=\psi(x+y+1)+1/(x+y+1)$, we have 
\begin{align*}
\frac{\partial g(x,y)}{\partial y} &=  \frac{\psi(x+y+2)-\psi(y+1)}{x+1}-\frac{1}{2(x+y+1)} 
-\frac{\psi(x+y+1)-\psi(y+1)}{x}+\frac{1}{2(x+y)}\\
 &=\frac{1}{x+y+1}\left[\frac{1}{x+1}+\frac{1}{2(x+y)}\right]
 -\frac{\psi(x+y+1)-\psi(y+1)}{x(x+1)}\\
 &\equiv\frac{h(x, y)}{x(x+1)}
\end{align*}
where 
$$h(x, y)=\frac{1}{x+y+1}\left[x+\frac{x(x+1)}{2(x+y)}\right]-\psi(x+y+1)+\psi(y+1).$$
Lemma \ref{lem1} records a useful monotonicity property of $h(x, y)$.
\begin{lemma}
\label{lem1}
For $x\geq 1/\sqrt{2}$ and fixed $y>0$, the function $h(x, y)$ strictly decreases in $x$.
\end{lemma}
{\bf Proof.} We have
\begin{align*}
\frac{\partial h(x,y)}{\partial x} &= \frac{y+1}{(x+y+1)^2} +\frac{2xy(x+y+1)+y^2+y}{2(x+y)^2(x+y+1)^2} 
-\psi'(x+y+1)\\
&< \frac{y+1}{(x+y+1)^2} +\frac{2xy(x+y+1)+y^2+y}{2(x+y)^2(x+y+1)^2} -\frac{1}{x+y+1}-\frac{1}{2(x+y+1)^2}\\
&=\frac{(1-2x^2)y-2x^3-x^2}{2(x+y)^2(x+y+1)^2}
\end{align*}
where (\ref{bp}) is used in the inequality.  From the last expression it is obvious that $\partial h(x, y)/\partial x<0$ if 
$x\geq 1/\sqrt{2}$. \qed
\begin{lemma}
\label{h1}
If $1/\sqrt{2}\leq x<1$ then $h(x, y)>0;$ if $x>1$ then $h(x, y)<0.$
\end{lemma}

{\bf Proof.} This follows from Lemma \ref{lem1}, noting that $h(1, y)=0$. \qed

Now deal with the $0<x\leq 1/\sqrt{2}$ case.  Let
$$u(x, y)=h(x, y)/x.$$
We have
\begin{align*}
\frac{\partial u(x, y)}{\partial x} &= \frac{\psi(x+y+1)-\psi(y+1)}{x^2}
-\frac{\psi'(x+y+1)}{x}+\frac{y-1}{2(x+y)^2}-\frac{y+2}{2(x+y+1)^2}.
\end{align*}
Using the Taylor expansion 
$$\psi(y+1)=\psi(x+y+1)-x\psi'(x+y+1)+\frac{x^2}{2} \psi''(\xi)$$
with $y+1<\xi<x+y+1$, we obtain
\begin{align*}
\frac{\partial u(x, y)}{\partial x} &= 
-\frac{\psi''(\xi)}{2}+\frac{y-1}{2(x+y)^2}-\frac{y+2}{2(x+y+1)^2}\\
&< \frac{1}{2(y+1)^2}+\frac{1}{2(y+1)^3}+\frac{1}{4(y+1)^4} 
+\frac{y-1}{2(x+y)^2}-\frac{y+2}{2(x+y+1)^2}\\
&\equiv v(x, y)
\end{align*}
where (\ref{bpp}) is used in the inequality.

\begin{lemma}
\label{h3}
The function $v(x, y)$ increases in $x$ for all $x, y>0$.
\end{lemma}

{\bf Proof.} We have
\begin{align*}
\frac{\partial v(x, y)}{\partial x} &= 
\frac{y+2}{(x+y+1)^3}-\frac{y-1}{(x+y)^3}\\
&=\frac{3x(x+y)(x+y+1)+3x+2y+1}{(x+y)^3(x+y+1)^3}\\
&>0. \qed
\end{align*}

\begin{lemma}
\label{h2}
For $0<x\leq 1/\sqrt{2}$ we have $u(x, y)>0$ (and hence $h(x, y)>0$). 
\end{lemma}

{\bf Proof.}  We first show that $u(x, y)$ decreases in $x$ if $x\in (0, 
1/\sqrt{2}]$.  In view of Lemma \ref{h3}, we know
$$\frac{\partial u(x, y)}{\partial x} <v(x, y) \leq v(1/\sqrt{2}, y).$$
After laborious but straightforward calculations we get
\begin{align}
\label{v}
v(1/\sqrt{2}, y)
&=\frac{-2y^2-2(5-2\sqrt{2})y-(5-2\sqrt{2})}{16(y+1)^4(1/\sqrt{2}+y)^2(y+1+1/\sqrt{2})^2}\\
\nonumber
&<0,
\end{align}
noting that $y>0$.  Thus $u(x, y)$ decreases in $x$ 
when $x\leq 1/\sqrt{2}$.  We may use this and Lemma \ref{h1} to 
obtain
\begin{align*}
u(x, y) &\geq u(1/\sqrt{2}, y)\\
   &=\sqrt{2} h(1/\sqrt{2}, y)\\
   &> 0. \qed
\end{align*}

{\bf Remark 1.} Expression (\ref{v}) can also be verified by a symbolic computing package such as Maple or Mathematica.

\begin{lemma}
\label{glim}
We have
\begin{equation}
\label{lim}
\lim_{y\rightarrow\infty} g(x, y)=0.
\end{equation}
\end{lemma}
{\bf Proof.} For fixed $x>0$, the asymptotic formula (\ref{lgam}) gives 
\begin{align*}
f(x, y)&=\frac{(x+y+1/2)\log(x+y+1)-(y+1/2)\log(y+1)-x}{x}-\frac{1}{2}\log(x+y)+O(y^{-1})\\
&=\log(x+y+1)-1+\frac{y+1/2}{x} \log\left(1+\frac{x}{y+1}\right)-\frac{1}{2}\log(x+y)+O(y^{-1})\\
&=\log(x+y+1) -1+\frac{y+1/2}{y+1}-\frac{1}{2}\log(x+y)+O(y^{-1})\\
&=\frac{1}{2}\log(y)+O(y^{-1})
\end{align*}
as $y\rightarrow\infty$.  Since $g(x, y)=f(x+1, y)-f(x, y)$, it is clear that $g(x, y)\rightarrow 0$ as $y\rightarrow\infty$ 
for fixed $x>0$. \qed

\begin{lemma}
\label{g}
If $x>1$ then $g(x,y)>0$; if $0< x<1$ then $g(x, y)<0$.
\end{lemma}

{\bf Proof.} By Lemma \ref{h1}, when $x>1$ we have $h(x, y)<0$, which means that $g(x, 
y)$ strictly decreases in $y$.  By Lemma \ref{glim}, $g(x, y)>0$ for all $x>1$ and $y>0$.  Similarly, if $0<x<1$, then by 
Lemmas \ref{h1} and \ref{h2} we have $h(x, y)>0$, i.e., $g(x, y)$ strictly increases in $y$, which implies $g(x, y)<0$ in 
view of (\ref{lim}). \qed

Theorem \ref{main} then follows from Lemma \ref{g}.

{\bf Remark 2.} The statement of the problem by Chen and Qi (\cite{CQ}, Open Problem 1) does not explicitly mention the range 
$x>1$.  We emphasize that $x>1$ is necessary and sufficient for (\ref{eqn1}). 

\section*{Acknowledgments} The author thanks an anonymous reviewer for his/her valuable comments.


\begin{thebibliography}{10}
\bibitem{AS}
M. Abramowitz, I. A. Stegun (Eds.), {\it Handbook of Mathematical Functions}, Dover Publications, New York, 1964. 


\bibitem{A}
H. Alzer, On some inequalities for the gamma and psi functions, 
{\it Mathematics of Computation} {\bf 66} (1997) 373--389.

\bibitem{CQ}
C. P. Chen, F. Qi, Monotonicity results for the gamma function, {\it J. Inequal. Pure Appl. Math.} {\bf 4} (2003) Art. 
44.

\bibitem{GQ}
B.-N. Guo, F. Qi, Inequalities and monotonicity for the ratio of gamma functions, {\it Taiwanese
J. Math.} {\bf 7} (2003) 239--247.

\bibitem{K}
S. Koumandos, Remarks on some completely monotonic functions, {\it J. Math. Anal. Appl.} {\bf 324} (2006) 1458–-1461.

\bibitem{QCCG}
F. Qi, R.-Q. Cui, C.-P. Chen, B.-N. Guo, Some completely monotonic functions involving polygamma functions and
an application, {\it J. Math. Anal. Appl.} {\bf 310} (2005) 303–-308.

\end{thebibliography}
\end{document}